\documentclass{elsart}
\usepackage{amssymb,amsfonts,amsmath}

\newtheorem{theorem}{Theorem}[section]
\newtheorem{lemma}[theorem]{Lemma}
\newtheorem{proposition}[theorem]{Proposition}
\newtheorem{corollary}[theorem]{Corollary} 
\newtheorem{Definition}[theorem]{Definition}
\newtheorem{Example}[theorem]{Example}
\newenvironment{example}{\begin{Example}\em}%
	{\hfill$\Box$\smallskip\end{Example}}
\newenvironment{proof}{\noindent{\bf Proof.}}{\hfill$\Box$\medskip}
\newenvironment{oneshot}[1]{\noindent\begin{trivlist}%
	\item[\hskip\labelsep{\sc#1}]}%
	{\hfill$\Box$\smallskip\end{trivlist}\vspace*{-.6cm}}

\def\Maple{{\sc Maple}} 
 \def\NN{{\mathbb N}} \def\PP{{\mathbb P}}
  \def\ZZ{{\mathbb Z}}

\begin{document}
\begin{frontmatter}

\title{Multihomogeneous resultant formulae by means of complexes}

\author{Alicia Dickenstein\thanksref{AD}}
\thanks[AD]{Partially supported by Action A00E02 of the ECOS-SeTCIP French-Argentina
 bilateral collaboration, UBACYT  X052 and ANPCYT 03-6568, Argentina}
\address{Departamento de Matem\a'{a}tica, F.C.E y N.,
 UBA (1428) Buenos Aires, Argentina}
\ead{alidick@dm.uba.ar}

\author{Ioannis Z.~Emiris\thanksref{IE}}
\thanks[IE]{Partially supported by Action A00E02 of the ECOS-SeTCIP French-Argentina
 bilateral collaboration and FET Open European Project IST-2001-35512 (GAIA-II).}
\address{Department of Informatics \& Telecommunications,
 University of Athens, Greece 15771, and INRIA Sophia-Antipolis, France}
\ead{emiris@di.uoa.gr}

\begin{abstract}
The first step in the generalization of the classical theory of
homogeneous equations to the case of arbitrary support is to
consider algebraic systems with multihomogeneous structure.
We propose constructive methods for resultant matrices
in the entire spectrum of resultant formulae, ranging from pure
Sylvester to pure B\'ezout types, and including matrices of hybrid
type of these two.
Our approach makes heavy use of the combinatorics of multihomogeneous
systems, inspired by and generalizing certain
joint results by Zelevinsky, and Sturmfels or Weyman
\cite{StZe}, \cite{WeZe}.
One contribution is to provide conditions and algorithmic
tools so as to classify and construct the smallest possible
determinantal formulae for multihomogeneous resultants.
Whenever such formulae exist, we specify the underlying complexes
so as to make the resultant matrix explicit.
We also examine the smallest Sylvester-type matrices,
generically of full rank, which yield a multiple of the resultant.
The last contribution is to characterize the systems that admit a
purely B\'ezout-type matrix and show a bijection of such
matrices with the permutations of the variable groups.
Interestingly, it is the same class of systems admitting an optimal
Sylvester-type formula.
We conclude with examples showing the kinds of matrices that
may be encountered, and illustrations of our \Maple\  implementation.
\end{abstract}

\end{frontmatter}

\section{Introduction}\label{Sintro}

Resultants provide efficient ways for studying and solving polynomial
systems by means of their matrices.
This paper considers the sparse (or toric) resultant, which
exploits {\em a priori} knowledge on the support of the equations.
We concentrate on unmixed (i.e.\  with identical supports)
systems where the variables can be partitioned
into groups so that every polynomial is homogeneous in each group.
Such polynomials, and the resulting systems, are called {\em
multihomogeneous}.
Multihomogeneous structure is a first step away from the
classical theory of homogeneous systems towards fully
exploiting arbitrary sparse structure.
Multihomogeneous systems are encountered in several areas
including geometric modeling
(e.g.\  \cite{ChGoZh98}, \cite{Saxe97}, \cite{Zhan00}),
game theory and computational economics \cite{McMc97}.

Known sparse resultant matrices are  of different types.
On the one end of the spectrum are the {\em pure Sylvester-type}
matrices,
where the polynomial coefficients fill in the nonzero entries of the
matrix;
such is the coefficient matrix of linear systems, Sylvester's matrix
for univariate polynomials, and Macaulay's matrix for homogeneous
systems.
On the other end are the {\em pure B\'ezout-type} matrices,
i.e.\ matrices
where the coefficients of the {\em Bezoutian} associ\-ated to the input
polynomials
fill in the nonzero entries of the matrix, whereas
hybrid matrices, such as Dixon's, contain blocks of both pure types.
The examples in Section~\ref{Sexample} show the intricacy of
such matrices.
Hence the interest to describe them in advance in terms of
combinatorial
data, which allows for a structured matrix representation, based on
quasi-Toeplitz or quasi-Hankel structure~\cite{EmPa02snap},\cite{MouPan00}.

Our work builds on \cite{WeZe} and their study of multihomogeneous
systems through the determinant of a resultant complex.
First, we give  precise {\em degree vectors} together with
algorithmic methods for identifying and constructing
{\em determinantal formulae} for the sparse resultant, i.e.\
matrices whose determinant equals the sparse resultant.
The underlying resultant complex is made explicit and computational
tools are derived in order to produce the smallest such formula.
Second, we describe and construct the
smallest possible pure Sylvester matrices, thus generalizing the
results of \cite{StZe} and  \cite[Sec.~13.2, Prop.~2.2]{GKZ}, already
present in the interesting paper \cite{mcc33} pointed by one of
the referees.
The corresponding systems  include
all systems for which exact Sylvester-type matrices are known.
We consider more general Sylvester-type matrices, and show that
in the search of small formulae, these more general matrices are
not crucial.
The third contribution of this paper is to offer
sufficient and necessary conditions for systems to admit
purely B\'ezout determinantal formulae, thus generalizing
a result from \cite{ChtKap00}.
It turns out that these are precisely the same systems
admitting optimal Sylvester-type formulae, and this is nothing
but a special case of complexes with only two non-vanishing
cohomologies.
We also show a bijection of such matrices with the permutations on
$\{1,\dots,r\}$, where $r$ stands for the number of the variable
groups.
While constructing explicit B\'ezout-type formulae, we derive a
precise description of the support of the Bezoutian polynomial.

The complex with terms $K_\nu(m)$ described in the next section
is known as the {\em Weyman complex}.
For any choice of  dimensions, of degrees of the input equations 
and of an integer vector $m$,
the multihomogeneous resultant equals the determinant of the 
Weyman complex (for the
corresponding
monomial basis at each of the terms), which can
be expressed as a quotient of products of subdeterminants
extracted from the differentials in the complex.
This way of defining the resultant was introduced by Cayley
\cite[app.~A]{GKZ}, \cite{Weym94}, \cite{WeZe}. In the particular
case in which the complex has just two terms, its determinant
is nothing but the determinant of the only nonzero differential,
which is therefore equal to the resultant. In this case, we say
that there is a {\it determinantal\/} formula for the resultant and the
corresponding
degree vector $m $ is called {\it determinantal\/}.
In \cite{WeZe}, the multihomogeneous systems for which a determinantal
formula exists were classified; see also \cite[Sec.~13.2]{GKZ}.
Their work, though, does not identify completely the corresponding
morphisms nor the determinantal vectors $m$,
a question we partially undertake. We follow the results in
\cite{DADi01}, which concerned the homogeneous case, inspired also
by \cite{Joua97}.

The main result of \cite{StZe}  was to prove that a determinantal
formula of Sylvester type  exists exactly when all defects are zero.
In \cite[Thm.~2]{StZe} (recalled in \cite[Sec.~13.2, Prop.~2.2]{GKZ};
see also \cite[Thm.~4]{mcc33}) all such formulae are characterized
by showing a bijection with the permutations of $\{1,\dots,r\}$ and
defined
the corresponding degree vector $m$ as in Definition \ref{Dmpi} below.
This includes all known Sylvester-type formulae, in particular,
linear systems, systems of 2 univariate polynomials and bihomogeneous
systems of 3 polynomials whose resultant is, respectively, the
coefficient determinant, the Sylvester resultant and the Dixon
resultant.
In fact, Sturmfels and Zelevinsky characterized all
determinantal Cayley-Koszul complexes, which are instances of the
Weyman complexes when all the higher cohomologies vanish.

The incremental algorithm for
sparse resultant matrices \cite{EmCa95} relies
on the determination of a degree vector $m$.
When $\delta=0$, it produces
optimal Sylvester matrices by \cite{StZe}.
For other multihomogeneous systems, \cite{EmCa95}
heuristically produces small matrices, yet with no guarantee.
For instance, on the system of Ex.~\ref{Ex211_222},
it finds a $1120\times 1200$ matrix.
The present paper explains the behavior of the algorithm,
since the latter uses degree vectors following
(\ref{Dmpi}) defined by random permutations.
Our results provide immediately
the smallest possible matrix.
More importantly, the same software constructs
all Sylvester-type formulae described here.

Pure B\'ezout-type formulae were studied in \cite{ChtKap00}
for unmixed systems whose support is the direct sum of what they call
basis simplices, i.e.\  the convex hull of the origin and another $l_k$
points, each lying on a coordinate axis.
This includes the case of multihomogeneous systems.
They showed that in this case a sparse resultant matrix
can be constructed from the Bezoutian polynomial, to be
defined in Section~\ref{Sbezout},
though the corresponding matrix formula is not always determinantal.
Their Corollary~4.2.1 states that for multihomogeneous systems with
null
defect vector, the B\'ezout formula becomes determinantal;
Saxena had proved the special case of all $l_k=1$ \cite{Saxe97}.
In \cite[Sec.~4.2]{ChtKap00} they indicate there are $r!$ such
formulae and in Section~5 they study bivariate systems ($n=2$)
showing that then, these are the only determinantal formulae.

Section~\ref{Sbezout} proves these results in a different manner
and characterizes the determinantal cases for
multihomogeneous systems, showing that a null defect vector is a
sufficient
but also {\em necessary} condition for a determinantal formula of pure
B\'ezout type for any
$n$. Thus, there is an optimal Sylvester-type formula for the resultant
if and only of there is an optimal pure B\'ezout-type formula.
 (cf.\ Definition \ref{defbezfla})
This had been proven for arbitrary systems only in the bivariate case
\cite{ChGoZh98}.
In particular, we explicitly exhibit a choice of the differential in
the Weyman complex in this case (cf.\ Theorem \ref{explicitbez}), thus
partially answering the ``challenge to make these maps explicit" of
\cite[Sec.~5.1]{WeZe}. 

Studies exist \cite{ChGoZh98}, \cite{DADi01}, \cite{Zhan00}
for dealing with
hybrid formulae including B\'ezout-type blocks or pure B\'ezout
matrices,
and concentrate on the computation of such matrices.
In particular, \cite{Zhan00} elaborates on the relation of Sylvester
and B\'ezout-type matrices (called Cayley-type there) and the
transformations that link them. The theoretical setting together
with Pfaffian formulae for resultants is addressed in \cite{EisSch01}.
This is made explicit for any toric surface in \cite{Kh02}.
In the recent preprint \cite{acg02}, not also the multihomogeneous resultant
but the whole ideal of inertia forms is studied, extending results
of Jouanolou \cite{j80} in the homogeneous case.

This paper is organized as follows.
The next section provides some technical facts useful later.
Section~\ref{Sdet} offers bounds in searching for the smallest
possible determinantal (hybrid) formulae.
Section~\ref{degreevectors} makes explicit one degree vector
attached to any determinantal data of dimensions of projective
spaces, and discusses further techniques for obtaining determinantal
formulae.
Section~\ref{Ssylvester} and \ref{Sbezout} characterize
matrices of pure Sylvester and pure B\'ezout type respectively.
In Section~\ref{Sexamples} we fully describe the formulae of a
system of 3 bilinear polynomials.
Then, we provide an explicit example of
a hybrid resultant matrix for a multidegree for which neither
pure Sylvester nor pure B\'ezout determinantal formulae exist;
this example illustrates the possible morphisms that may be
encountered with multihomogeneous systems.
Our \Maple\  implementation is described in Section~\ref{Scode}.

A preliminary version of certain results in this paper has
appeared in \cite{DicEmi02}.

\section{Preliminary observations} \label{Sprelim}

We consider the $r$-fold product
$X:=\PP^{l_1} \times \ldots \times \PP^{l_r}$ of projective
spaces of respective dimensions $l_1,\dots,l_r$ over an
algebraically closed field of characteristic zero,
for some natural number $r$.
We denote  by $n=\sum_{k=1}^r l_k$ the dimension of $X$, i.e.\
the number of affine variables.

\begin{Definition}
Consider $d=(d_1,\ldots,d_r)\in\NN_{>0}^r$ and multihomogeneous
polynomials $f_0,\ldots, f_n$ of degree $d$.
The multihomogeneous resultant
is an irreducible polynomial
$R(f_0,\ldots,f_n) = R_{(l_1,\ldots,l_r),d}(f_0,\ldots,f_n)$
in the coefficients of
$f_0,\ldots,f_n$ which vanishes iff the
polynomials have a common root in $X$.
\end{Definition}

This is an instance of the sparse resultant~\cite{GKZ}.
It may be chosen with integer coefficients, and it is
uniquely defined up to sign by the requirement that it has
relatively prime coefficients.
The resultant polynomial is itself homogeneous in the coefficients
of each $f_i$, with degree given by the multihomogeneous B\'ezout
bound ${\binom n {l_1,\dots,l_r}} d_1^{l_1}\cdots d_r^{l_r}$
\cite[Prop.~13.2.1]{GKZ}.
This number is also called the $m$-homogeneous bound
\cite{Wamp92}.

Let $V$ be the space of $(n+1)$ tuples $f=(f_0,\dots,f_n)$ of
multihomogeneous forms of degree $d$ over $X$. Given a
{\em degree vector} $m \in \ZZ^r$
 there exists a finite complex $K_\cdot = K_\cdot (m)$
of free modules over the ring of polynomial functions on
$V$ \cite{WeZe}, whose terms depend only on
$(l_1,\ldots,l_r), d$ and $m$ and whose differentials
are polynomials on $V$ satisfying:\\
(i) For every given $f$ we can specialize the
differentials in $K_\cdot$ by evaluating at $f$ to get
a complex of finite-dimensional vector spaces.\\
(ii) This complex is exact iff
$R(f_0,\ldots,f_n) \not= 0$.\\

In order to describe the terms in these complexes
some facts from cohomology theory are necessary;
see \cite{Hart77} for details.
Given a {\em degree vector} $m \in \ZZ^r$,
define, for $\nu\in\{-n,\dots,n+1\}$,
\begin{equation} \label{EKnu_coH}
K_{\nu}(m) =  \bigoplus_{p\in \{0,\ldots,n+1\}}
 H^{p-\nu}(X, m-pd)^{{\binom{n+1} p}},
\end{equation}
where for an integer $r$-tuple $m'$, $H^q(X, m')$
denotes the $q$-th cohomology of $X$ with coefficients
in the sheaf ${\mathcal O}(m')$ such that its global sections
$H^0(X, m')$ are identified with  multihomogeneous  polynomials of
(multi)degree $m'$.
By the K\"unneth formula, we have
$$
H^q(X, m-pd) =  \bigoplus_{j_1+\cdots+j_r=q}^{j_k\in\{0,l_k\}}
 \bigotimes_{k=1}^r H^{j_k}(\PP^{l_k}, m_k-pd_k),
$$
where $q=p-\nu$ and the second sum runs over all integer sums
$j_1+\cdots+j_r=q, j_k\in\{0,l_k\}$. In particular,
$H^0(\PP^{l_k}, \alpha_k)$ is the space of all homogeneous polynomials
in $l_k+1$ variables with total degree $\alpha_k$.  By Serre's duality,
for any $\alpha \in \ZZ^r$, we also know that
\begin{equation} \label{Serre}
H^q (X, \alpha)  \simeq H^{n-q}(X, (-l_1-1, \dots, -l_r-1) - \alpha)^*,
\end{equation}
where ${}^*$ denotes dual.
We recall Bott's formulae for these cohomologies.
\begin{proposition}
For any $m\in\ZZ^r$,
$H^{l_k}(\PP^{l_k}, m_k-pd_k)=0\Leftrightarrow m_k-pd_k\ge -l_k$,
$H^{0}(\PP^{l_k}, m_k-pd_k)=0\Leftrightarrow m_k-pd_k<0$, for
$k\in\{1,\dots,r\}$. Moreover,
\begin{eqnarray*}
H^{j}(\PP^{l_k}, m_k-pd_k) =0 , \forall j \not= 0, l_k, \\
\dim H^{l_k}(\PP^{l_k}, m_k-pd_k) = {\binom{-m_k+pd_k-1}{l_k}},\\
\dim H^{0}(\PP^{l_k}, m_k-pd_k) = {\binom{m_k-pd_k+l_k}{l_k}}.\\
\end{eqnarray*}
\end{proposition}
Consequently,
$$
\dim H^q(X, m-pd) = \sum_{j_1+\cdots+j_r=q}^{j_k\in\{0,l_k\}}
 \prod_{k=1}^r \dim H^{j_k}(\PP^{l_k}, m_k-pd_k),\\
$$
and
$$
\dim K_{\nu}(m) = \sum_{p\in [0,n+1]}
 {\binom{n+1} p} \dim H^{p-\nu}(X, m-pd).
$$

\begin{Definition}\label{Dcrit}
Given $r$ and $(l_1,\dots,l_r), (d_1,\dots,d_r)\in\NN^r$,
define the {\em defect} vector $\delta\in\ZZ^r$
(just as in \cite{StZe}, \cite{WeZe})
by $\delta_k:=l_k - \lceil \frac{l_k}{d_k} \rceil$.
Clearly, this is a non-negative vector.
We also define the {\em critical degree} vector
$\rho\in\NN^r$ by $\rho_k:=(n+1)d_k -l_k -1,$ for all $k=1,\ldots r.$
\end{Definition}

\begin{lemma} {\rm \protect\cite{WeZe}} \label{Ldefect}
For any $i\in [r] := \{1,\dots,r\},\,
d_i l_i<d_i+l_i \Leftrightarrow \delta_i=0 \Leftrightarrow
\min\{l_i,d_i\}=1$.
\end{lemma}

Let us establish a general technical lemma.

\begin{lemma} \label{Ldeltas}
For any $k\in \{1,\dots,r\},\, 0\leq (l_k-\delta_k)d_k-l_k\leq d_k-1$.
\end{lemma}

\begin{proof}
By definition, $\delta_k=l_k-\lceil l_k/d_k\rceil \Leftrightarrow
\lceil l_k/d_k\rceil=l_k-\delta_k = (l_k+t_k)/d_k$ for some integer
$t_k$ such that $0\leq t_k\leq d_k-1$.
\end{proof} 

We detail now the main results in \cite{WeZe}.
They show (Lemma 3.3(a)) that
a vector $m\in\ZZ^r$ is determinantal iff
$K_{-1}(m)=K_2(m)=0$.  They also prove  in Theorem~3.1 that
a determinantal vector $m$ exists iff $\delta_k \le 2$ for all
$k \in [r]$. To describe a differential in the complex from $K_\nu(m)$ to
$K_{\nu +1}(m)$, one needs to describe all the morphisms
$\delta_{p,p'}$ from the summand
corresponding to an integer $p$ to the summand corresponding to
another integer $p'$, where both $p,p' \in\{0,\dots, n+1\}$.
\cite[Prop.~2.5, 2.6]{WeZe} proves this map is $0$ when $p<p'$
and that, roughly speaking, it
corresponds to a {\em Sylvester} map $(g_0,\ldots, g_n)  \to
\sum_{i=0}^n g_i f_i$ when $p= p' +1$, thus having all nonzero
entries in the corresponding matrix given by coefficients of
$f_0,\dots, f_n$. For $p > p' +1$, the maps $\delta_{p,p'}$
are called higher-order differentials. By degree reasons, they cannot be
given
by Sylvester matrices. Theorem~2.10 also gives
an explicit theoretical construction
of the higher-order differentials in the pure B\'ezout case
(cf.\ Definition \ref{defbezfla}).

\section{Bounds for determinantal degree vectors} \label{Sdet}

This section addresses the computational problem of
enumerating all determinantal degree vectors $m\in\ZZ^r$.
The ``procedure" of \cite[Sec.~3]{WeZe} ``is quite
explicit but it seems that there is no nice way to
parametrize these vectors", as stated in that paper.
Instead, we bound the range of $m$ to
implement a computer search for them. In Section~\ref{degreevectors}
we will give an explicit choice of degree vector $m$ for each
determinantal data $(~l_1,\dots,l_r; d_1,\dots,d_r)$.

Given $k\in\{1,\dots,r\}$ and a vector $m \in \ZZ^r$ define as in
\cite{WeZe}:
$$
P_k(m)= \left\{p\in\ZZ : \frac{m_k}{d_k} < p\le \frac{m_k+l_k}{d_k}
\right\}.
$$
Let $\tilde{P}_k(m)$ be the real interval
$\left(\frac{m_k}{d_k},\frac{m_k+l_k}{d_k}\right]$,
so $P_k(m)=\tilde{P}_k(m)\cap\ZZ$.
Using Lemma~3.3 in \cite{WeZe}, it is easy to give bounds for all
determinantal vectors $m$ for which all $P_k(m)\neq\emptyset$.
\begin{lemma} \label{Lboundmk}
For a determinantal $m\in\ZZ^r$ and
for all $k\in\{1,\dots,r\}, \, P_k(m) \neq\emptyset$ implies,
$\max\{-d_k,-l_k\} \le m_k \le d_k(n+1) -1 + \min\{d_k-l_k,0\}.$
\end{lemma}

\begin{proof} 
By \cite[Lem.~3.3(b)]{WeZe}, $P_k(m) \subset[0,n+1] \Rightarrow
m_k/d_k\ge -1$ and $(m_k+l_k)/d_k\ge 0$, which imply the lower bound.
Also, $(m_k+l_k)/d_k<n+2\Leftrightarrow m_k\le (n+2)d_k-l_k-1$
and $m_k/d_k<n+1\Leftrightarrow m_k\le (n+1)d_k-1$ yield the upper
bound.
Notice that the possible values for $m_k$ form a non-empty set,
since the two bounds are negative and positive respectively.
\end{proof}

Now, $p\in P_k(m)$ iff
$H^{l_k}(\PP^{l_k}, m_k-pd_k) = H^{0}(\PP^{l_k}, m_k-pd_k)=0.$
Thus, a first guess could  be that  all determinantal
vectors give $P_k(m)\neq \emptyset$. But this
is not the case, as the following example shows:

\begin{example} \label{Ex12_23} 
Set $l=(1,2), d=(2,3)$.
We focus on degree vectors of the form
$m=(2\mu_1,3\mu_2)$, for $\mu_1,\mu_2\in\ZZ$.
Then, for all such $m$, the sets $P_k(m)$, for $k=1,2$,
are empty. Nevertheless, there exist  $4$ determinantal
vectors of this form, namely $m=(4,3), (0,6), (2,6)$
or $(6,3)$.  Moreover, the vector $m =(6,3)$
gives a determinantal formula with a matrix of size $88$,
which is closer to the smallest possible one which has size
$72$. The largest determinantal formula is given
by the determinant of a square matrix of size $180$.
Note that the degree of the multihomogeneous resultant
is $216$.
More details on this example are provided in Section \ref{Scode}.
\end{example}

We wish now to get a bound for those determinantal vectors
for which some $P_k(m)$ is empty.
Let $[\cdot]_k \in\{0,1,\dots,d_k-1\}$
denote remainder after division by $d_k$.

\begin{Definition} \label{Dpertm}
Given $m\in\ZZ^r$ and $k \in \{1,\ldots,r\}$, define new vectors
$m',m''\in\ZZ^r$ whose $j$-th coordinates equal those of $m$
 for all $j \neq k$ and  such that
$m_k'= m_k + d_k - [m_k]_k -1 \geq m_k \geq m_k''=m_k - [m_k+l_k]_k.$
\end{Definition}

\begin{lemma} \label{Lwelldef}
The  vectors $m',m''$ differ from $m$ at their $k$-th
coordinate if $P_k(m)=\emptyset$.
\end{lemma}

\begin{proof} 
Let us write $m_k=jd_k+[m_k]_k$ for $j\in\ZZ$. Then,
$\frac{m_k}{d_k}=j+\frac{[m_k]_k}{d_k}$. If
$P_k(m) =\emptyset, \,
\frac{m_k+l_k}{d_k}<j+1\Rightarrow jd_k+[m_k]_k+l_k<(j+1)d_k$,
so $[m_k]_k\le d_k-l_k-1 \Rightarrow l_k \le -[m_k]_k+d_k-1$
and thus $1\le d_k - [m_k]_k-1$.
Also, $[m_k+l_k]_k\ge 1$
because $[m_k+l_k]_k=0\Rightarrow (m_k+l_k)/d_k \in P_k(m)$.
\end{proof}

\begin{lemma} \label{Lpertmk1}
If $m\in\ZZ^r$ with $P_k(m)=\emptyset$ and
$H^0(\PP^{l_k}, m_k-pd_k)=0$  (resp. $H^{l_k}(\PP^{l_k}, m_k-pd_k)=0$),
then $P_k(m')\neq\emptyset$ and $H^0(\PP^{l_k}, m_k'-pd_k)=0$ (resp.
$H^{l_k}(\PP^{l_k}, m_k'-pd_k)=0$), where $m_k'=m_k + d_k- [m_k]_k-1$ as
in Definition \ref{Dpertm}.
\end{lemma}

\begin{proof} 
Write $m_k=jd_k +[m_k]_k$ for some integer $j\in \ZZ$.
To prove $P_k(m')\neq\emptyset$ we show $j+1\in P_k(m')$, i.e.\
$\frac{m_k'}{d_k}< j+1 \le \frac{m_k'+l_k}{d_k}$ $\Leftrightarrow
m_k'<(j+1)d_k\le m_k'+l_k$  $\Leftrightarrow$
$jd_k + d_k-1 < (j+1)d_k \le jd_k + d_k-1+l_k, $
which is clearly true since $1\le l_k$.
Now, $H^{l_k}(\PP^{l_k}, m_k-pd_k)=0\Leftrightarrow m_k+l_k\ge pd_k$
hence $m_k'+l_k\ge pd_k$ because $m_k'\ge m_k$.
$H^0(\PP^{l_k}, m_k-pd_k)=0\Leftrightarrow m_k<pd_k$ so $j\le p-1$.
Hence, $m_k-[m_k]_k=jd_k\le (p-1)d_k\Leftrightarrow m_k-[m_k]_k+d_k\le
pd_k$
which is the desired conclusion.
By \cite{WeZe}, $P_k(m')\subset [0,n+1]$ from which
$j\in\{-1,0,\dots,l\}$.
\end{proof}

\begin{lemma} \label{Lpertmk2}
If $m\in\ZZ^r$ with $P_k(m)=\emptyset$ and
$H^0(\PP^{l_k}, m_k-pd_k)=0$ (resp. $H^{l_k}(\PP^{l_k}, m_k-pd_k)=0$),
then $P_k(m'')\neq\emptyset$ and $H^0(\PP^{l_k}, m_k''-pd_k)=0$
(resp. $H^{l_k}(\PP^{l_k}, m_k''-pd_k)=0$), where $m_k''=m_k -
[m_k+l_k]_k$ as in Definition \ref{Dpertm}.
\end{lemma}

\begin{proof} 
Write $m_k+l_k=jd_k+[m_k+l_k]_k$ for some integer $j\ge 0$.
To prove $P_k(m_k'')\neq\emptyset$ we show it contains
$j$, i.e.,
$\frac{m_k - [m_k+l_k]_k}{d_k}< j \le \frac{m_k - [m_k+l_k]_k+l_k}{d_k}$ $\Leftrightarrow$
$m_k - [m_k+l_k]_k<jd_k\le m_k - [m_k+l_k]_k+l_k$ $\Leftrightarrow$ $jd_k - l_k < jd_k \le jd_k,$
which is clearly true since $0< l_k$.
$H^0(\PP^{l_k}, m_k-pd_k)=0\Leftrightarrow m_k<pd_k\Rightarrow m_k'' < pd_k$
because $m_k''\le m_k$, hence $H^0(\PP^{l_k}, m_k''-pd_k)=0$.
$H^{l_k}(\PP^{l_k}, m_k-pd_k)=0\Leftrightarrow m_k+l_k\ge pd_k$, then $j\ge p$.
Hence $m_k+l_k-[m_k+l_k]_k\ge p d_k$ which finishes the proof.
By \cite{WeZe}, $P_k(m'')\subset [0,n+1]$ hence $j\in\{0,\dots,n+1\}$.
\end{proof}

Lemmas~\ref{Lpertmk1} and  \ref{Lpertmk2} imply

\begin{theorem}
For any determinantal $m\in\ZZ^r$, define vectors $m',m''\in\ZZ^r$
as in Definition \ref{Dpertm}
which differ from $m$ only at the $k$-th coordinates, $1\le k\le r$,
such that $P_k(m)=\emptyset$.
Then $P_k(m')\neq\emptyset, P_k(m'')\neq\emptyset$ and both
$m',m''$ are determinantal.
\end{theorem}

\begin{corollary} \label{Cbnd_empt}
For a determinantal $m\in\ZZ^r$ with $P_k(m)=\emptyset$
for some $k\in\{1,\dots,r\}$, we have
$0 \le m_k \le d_k(n+1) - l_k -1.$
\end{corollary}

\begin{proof}
Since $m_k',m_k''$ define $P_k(m')\neq\emptyset, P_k(m'')\neq\emptyset$,
we can apply Lemma~\ref{Lboundmk}.
We use the lower bound with $m_k''$ because $m_k'' <m_k<m_k'$.
$P_k(m)=\emptyset \Rightarrow d_k>l_k$, so
$m_k''=m_k-[m_k+l_k]_k \ge -l_k
\;$  $ \Rightarrow$  $\; m_k\ge [m_k+l_k]_k-l_k\ge 1-l_k,$
because $[m_k+l_k]_k\ge 1$ by the proof of Lemma~\ref{Lwelldef}.
If $m_k<0$, for $P_k(m)$ to be empty we need
$m_k+l_k<0 \Leftrightarrow m_k<-l_k$ which contradicts the derived
lower bound; so $m_k\ge 0$.
For the upper bound,
$m_k'=m_k+d_k-[m_k]_k-1 \le d_k(n+1)-1 \;$  $\Rightarrow$
$m_k\le d_k n+[m_k]_k \le d_k(n+1)-l_k-1;$
the latter follows from $[m_k]_k < d_k-l_k$,
(Lemma~\ref{Lwelldef}).
$(m_k+l_k)/d_k \le n+1-(1/d_k)<n+1$
implies the inclusion of the half-open interval in $(0,n+1)$.
The possible values for $m_k$ form a non-empty set,
since the lower bound is zero and the upper bound is
$d_k(n+1) - l_k -1\ge d_k-1>0$ since $d_k>l_k\ge 1$.
\end{proof}
So, in fact, the real interval $\tilde{P}_k \subset (0, n+1).$

\begin{corollary}\label{detbounds}
For a determinantal $m\in\ZZ^r$ and  $ k\in [r]$,
$\max\{-d_k,-l_k\} \le m_k \le d_k(n+1) - 1 + \min\{d_k-l_k,0\}.$
\end{corollary}

This implies there is a finite number of vectors
to be tested in order to enumerate all possible determinantal $m$.
This could also be deduced from the fact that the dimension of $K_0(m)$
equals the degree of the resultant.
Corollary~\ref{detbounds} gives a precise bound for the box in which to
search algorithmically for all determinantal $m$, including
those that are ``pure" in the terminology of \cite{WeZe}.
Our \Maple\  implementation, along with examples,
is presented in Section \ref{Scode}.

If we take $r=1=l_1$ and $m=2d_1-1$ we obtain the classical
Sylvester formula and the upper bound
given by Corollary \ref{Cbnd_empt} is attained.
At the lower bound $m=-1$, the corresponding complex
$H^1(\PP^{1},-1-2d)=H^0(\PP^{1}, 2d_1 -1)^* \to H^1(\PP^{1}, d_1-1)^2 = (H^0(\PP^{1}, d_1-1)^*)^2$
yields the same matrix transposed.
In addition, the system of 3 bilinear polynomials in
Section \ref{Sbilinear} admits a pure Sylvester formulae with
$m=(2,-1)$, which attains both lower and upper bounds of the Corollary.
The bounds in Lemma \ref{Lboundmk} can also be attained
(see Ex.\ \ref{Ex211_222} continued in Section \ref{Scode})
hence Corollary~\ref{detbounds} is tight.
It is possible that some combination of the coordinates
of $m$ restricts the search space.

\section{Explicit determinantal degree vectors}  \label{degreevectors}

We focus on the case of $\delta_k\le 2$ for all  $k$, which is
a necessary and sufficient condition for the existence
of a determinantal complex \cite{WeZe}.
We shall describe specific degree vectors $m$ that yield
determinantal complexes.

\begin{theorem} \label{Tdetmvec}
Suppose that $0\le \delta_k\le 2 $ for all $k=1, \dots, r$
and $\pi:[r]\to [r]$ is any permutation.
Then, the degree vector $m^{\pi} \in\ZZ^r$ with
\begin{equation} \label{syldetvec}
m_k^{\pi} = \left(1-\delta_k +
\sum_{\pi(j)\ge\pi(k)} l_j \right) d_k -l_k,
\end{equation}
for $k=1,\dots,r$ defines a determinantal complex.
\end{theorem}

\begin{proof}
First, $K_1\neq 0$ and $K_{0}\neq 0$ because they each contain
at least one nonzero direct summand, namely $H^0(X, m^{\pi}-pd)$
for $p=1,0$ respectively.
To see this, it suffices to prove $m^{\pi}_k-d_k\ge 0,\, k=1,\dots,r$,
which follows from $(1-\delta_k +l_k ) d_k -l_k-d_k\ge 0
\Leftrightarrow (-\delta_k +l_k ) d_k -l_k\ge 0$.
This holds by Lemma \ref{Ldeltas}.

To demonstrate that $K_2=0$ we shall see that every one
of its direct summands $H^q(X, m^{\pi}-(q+2)d)$ vanishes for
$q=\sum_{\pi(j)\in J} l_j$, where $J$ is any proper subset of $[r]$.
Recall that the case $J=[r]$ is irrelevant by the definitions of
Section \ref{Sprelim}.
It is enough to show $H^{0}(\PP^{l_k}, m^{\pi}_k-(q+2)d_k)=0$,
for some $k$ with $\pi(k)\not\in J$.
Let  $k$ be the maximum of the indices verifying $\pi(k)\not\in J$;
such a $k\in [r]$ always exists because $J\neq [r]$.
Since $(1 -\delta_k+l_k )d_k-l_k < 2 d_k$ by Lemma
\ref{Ldeltas}, it easily follows that
\begin{equation} \label{EstepK2}
m^{\pi}_k < (2+\sum_{\pi(j)\in J} l_j)d_k.
\end{equation}
It now suffices to establish $K_{-1}(m)=0$.
Since this module has no zero cohomology summand,
let $q=\sum_{\pi(j)\in J} l_j$, with
$J \not= \emptyset$. Let  $k\in J$ such that
$\pi(j) \ge \pi(k) $  for all $j \in J$.
Then, $\delta_k \le 2$ implies
$(2 -\delta_k)d_k \ge 0$, from which
$m_k^{\pi}-(q-1)d_k \ge -l_k ,$
and so $H^q(X, m^{\pi}-(q-1)d) =0$ for any direct summand of $K_{-1}(m)$.
\end{proof}

When the defects are at most $1$, we can give another explicit choice of
determinantal degree vector for each permutation of $[r]$.

\begin{theorem} \label{Tdetmvec2}
Suppose that $0\le \delta_k\le 1$ for all  $k=1, \dots, r$
and $\pi:[r]\to [r]$ is any permutation.
Then, the degree vector $m^{\pi} \in\ZZ^r$ with
\begin{equation} \label{bezdetvec}
m_k^{\pi} = \left(-\delta_k +
\sum_{\pi(j)\ge\pi(k)} l_j \right) d_k -l_k,
\end{equation}
for $k=1,\dots,r$ defines a determinantal complex.
\end{theorem}

\begin{proof}
Let us modify the previous proof.
First, $K_{0}\neq 0$ because it contains
$H^0(X, m^{\pi}-pd)\neq 0$ for $p=0$.
This follows from $(-\delta_k +l_k ) d_k -l_k\ge 0
\Leftrightarrow (-\delta_k +l_k ) d_k -l_k\ge 0$ for all $k$.
This holds by Lemma \ref{Ldeltas}.

Next, $K_1\neq 0$ because for $p=n+1$,
$H^{l_k}(\PP^{l_k}, m_k^{\pi}-pd_k)\neq 0$ for all $k$.
This follows from $(-\delta_k + n ) d_k -l_k < (n+1)d_k-l_k
\Leftrightarrow -\delta_k < 1$.

To demonstrate that $K_2=0$ we repeat the corresponding argument
in the proof of Theorem \ref{Tdetmvec}.
Then, it suffices to show that
$ (-\delta_k+l_k )d_k-l_k < 2 d_k, $
which is weaker than the inequality before (\ref{EstepK2}) above.
Similarly, to establish $K_{-1}(m)=0 \,$ the steps in the previous
proof lead us to showing
$ m_k^{\pi}-(q-1)d_k \ge -l_k \Leftarrow (1 -\delta_k)d_k \ge 0 $
which holds for $\delta_k \le 1$.
\end{proof}

One can easily verify that the vectors $m^\pi$ in both
Theorems above satisfy the bounds of Corollary~\ref{detbounds}.

A natural question is how to give an explicit  degree
vector yielding a determinantal formula of smallest size.  When $r=1$, this
is the content of Lemma 5.3 in \cite{DADi01}. Even for $r=2$, this seems
to be a difficult task in general.

\begin{example}
Set $l=(2,2), d=(3,2)$. Let $\pi_1:[2]\to[2]$ denote the identity and
$\pi_2:[2]\to[2]$ the permutation which interchanges $1$ and $2$.
The two determinantal vectors defined according to (\ref{syldetvec})
are $m_{\pi_1} = (10,2)$ and $m_{\pi_2} = (4,6)$, which yield determinantal
matrices of respective sizes $396$ and $420$.  The two
determinantal vectors defined according to (\ref{bezdetvec})
are $m_{\pi_1} = (7,0)$ and $m_{\pi_2} = (1,4)$, which yield determinantal
matrices of respective sizes $756$ and $780$.  But as we will see
in \S \ref{Sbezout}, these latter formulae give instead the vectors
providing
the smallest determinantal formulas when all defects are $0$.
These vectors can be computed by the function {\tt comp\_m} of
our implementation, discussed in Section~\ref{Scode}.

Also, when $r=1$ it is shown in \cite{DADi01} that the smallest formula
is attained for ``central'' determinantal degree vectors. In this
example,  all degree vectors $m$ satisfy that either $1 \leq m_1 \leq 5$
and $4 \leq m_2 \leq 7$ or $7 \leq m_1 \leq 11$ and $0 \leq m_2 \leq 4$;
these bounds are computed by the routines described in \S \ref{Scode}.
The smallest resultant matrix has size $340 \times 340$, and
corresponds to the degree vectors $(3,6)$ or $(9,1)$, which are in a
``central'' position among determinantal degree vectors, in other words,
their coordinates lie in a ``central'' position between the respective
coordinates of  other determinantal vectors.  
However, the  vectors $(k,4)$ are determinantal for $k$ from $1$ to $5$, but
the size of a matrix corresponding  to the  vector $(3,4)$ (both of whose
coordinates lie between the coordinates of the vectors $(1,4)$ and
$(5,4)$), equals $580$. So, it is bigger than the size of a resultant matrix
associated to the vector $(5,4)$, which is of dimension $540$.
\end{example}

In the above example, the two degree vectors giving the smallest
resultant matrices satisfy $(9,1) + (3,6) = (12,7)$, which is the
critical vector from Definition \ref{Dcrit}.
It is clear that for any determinantal $m$, the vector $(12,7) -m$
is also determinantal yielding the same matrix dimension.
This is a consequence of Serre's duality recalled in (\ref{Serre}).
The general statement in summarized in the next Proposition.

\begin{proposition}
Assume $m, m'\in \ZZ^r$ satisfy $m + m'= \rho$, the latter being the
critical degree vector of Definition \ref{Dcrit}.
Then,  $K_\nu(m)$ is dual to $K_{1-\nu} (m')$ for all $\nu \in \ZZ.$  
In particular, $m$ is determinantal if and only if $m'$ is determinantal,
yielding matrices of the same size, namely $\dim(K_0(m)) = \dim (K_1(m')).$
\end{proposition}

\begin{proof}
Based on the equality $m +m'= \rho$ we deduce that for all $p =0, \ldots, n+1$,
it holds that $(m' -pd) = (-l_1-1, \dots, -l_r-1) - (m- (n+1-p)d).$ Therefore, 
for all $q=0,\dots,n$, Serre's duality (\ref{Serre}) implies that $H^q(X, m'-pd)$
and $H^{n-q}(X, m- (n+1-p)d)$ are dual. Since $(n+1-p) - (n-q) = 1 -(p-q)$, we deduce that
$K_\nu(m)$ is dual to $K_{1-\nu} (m')$ for all $\nu \in \ZZ$, as desired.
Observe that in particular $K_{-1}(m) \simeq K_2(m')^*$ and
$K_0(m) \simeq  K_1(m')^*$, the latter giving the matrix dimension 
in the case of determinantal formulae.
\end{proof}

We end this section making explicit a consequence of Proposition~3.7
in \cite{WeZe} (and giving an independent proof),
which gives a generalization of the characterization of determinantal
complexes
in \cite{StZe}, \cite{WeZe}.  This result will allow us to give in \S \ref{Sbezout}
explicit expressions for all degree
vectors yielding a determinantal formula of smallest size when all defects
vanish.

\begin{theorem} \label{twoterms}
There exists a determinantal vector $m$ such that the
Weyman complex is reduced to only one non-zero cohomology group
on each of $K_0(m), K_1(m)$ if and only if all defects vanish, i.e.\
$\delta_k=0$ for all $k=1,\dots,r$.
\end{theorem}

\begin{proof}
Recall that for each $p\in\{0,\dots,n+1\}$ there exists at most
one integer $j$ such that $H^{j}\left(X, m-pd \right) \not=0$ in $K_{\nu}(m)$
where $p=j-\nu$ \cite[Prop.~2.4]{WeZe}.
In fact, let  $A (p):= \{ k \, : \, m_k - p d_k < - l_k \}$
and $B(p):= \{k \, : \, m_k - p d_k \geq 0 \}$.
Denote $j(p) := \sum_{ k \in A(p)} l_k.$
Then $H^{j'}\left(X, m-pd \right) =0$
for all $j'\not=j(p)$ and $H^{j(p)}\left(X, m-pd \right) \not=0$
iff $A(p) \cup B(p) = \{1,\ldots,r\}.$

The assumption of the theorem means that there exist
exactly two integers $p_1, p_2 \in \{0,\dots,n+1\}$ for
which $H^{j(p_i)}\left(X, m-p_id \right) \not=0, i=1,2$;
cf.\ also \cite[Lem.~3.3(a)]{WeZe}.
Then, for any $p \in \{0,\ldots,n+1\} \setminus \{p_1,p_2\}$,
there exists
$ k$ such that $k \notin A(p) \cup B(p)$, i.e., $p \in P_k(m)$.
The latter set is defined at the beginning of Section \ref{Sdet}.
Then, $\{0,\ldots,n+1\} \setminus \{p_1,p_2\} \subseteq
\cup_{k=1}^r P_k(m)$ and so
$$
l_1 + \ldots + l_r = n \leq \sum_{k=1}^r \# P_k(m) \leq
\sum_{k=1}^r \lceil \frac{l_k}{d_k} \rceil,
$$
where the first inequality uses the fact that
$\# \cup_{k=1}^r P_k(m) \le \sum_{k=1}^r \# P_k(m)$
and the second follows from the definition of $P_k(m)$.
Since $l_k \geq \lceil \frac{l_k}{d_k} \rceil$ for all $k$, we deduce
that
 $l_k =\lceil \frac{l_k}{d_k} \rceil$, and this can only happen
iff $l_k=1$ or $d_k =1.$
\end{proof}

\section{Pure Sylvester-type formulae} \label{Ssylvester}

This section constructs rectangular matrices of pure Sylvester-type
that have at least one maximal minor which is a nontrivial multiple of
the sparse resultant, coming from a complex of the form:
$$
\dots \rightarrow K_2(m) \rightarrow K_1(m)
\rightarrow K_0(m) \rightarrow K_{-1}(m) =0,
$$
where $K_1(m)=H^j(X, m-d)^{n+1}$, and $K_0(m)= H^0 (X, m)$ for a non negative
vector $m \in \ZZ_{\geq 0}^r.$ 

We assume that $H^0(X, m-p d)\neq 0$ for $\nu=0,1$ and $p-\nu=0$.
This implies $H^0(\PP^{l_k}, m_k-\nu d_k)\neq 0$ for all
$k\in\{1,\dots,r\}$.
Moreover, we must have $H^{p-\nu}(X, m-pd)=0$
for $p- \nu =\sum_{j\in J} l_j$ where
$J$ is any subset satisfying $\emptyset\neq J\subset\{1,\dots,r\}$.
Note that we do not require in general that $K_2(m)\neq 0$.

\begin{lemma} \label{Lstairs}
If $m \not= m'\in\ZZ^r_{\geq 0}$ yield a Sylvester-type matrix and
$m_k' \ge m_k$ for all $k\in\{1,\dots,r\}$
then, the Sylvester matrix associated to $m'$  is strictly larger
than the Sylvester matrix associated to $m$.
\end{lemma}

\begin{proof} 
We must show $\dim K_0(m') \ge \dim K_0(m)$, i.e.\
$\dim H^{0}(\PP^{l_k}, m_k'-pd_k) \ge \dim H^{0}(\PP^{l_k}, m_k-pd_k)$
for $k\in\{1,\dots,r\}$, i.e.\
${\binom{m_k'-pd_k+l_k}{l_k}} \ge {\binom{m_k-pd_k+l_k}{l_k}}.  $
The cohomology is nonzero, thus $m_k'-pd_k\ge m_k-pd_k\ge 0$
and this implies the desired inequality because
${\binom{s+l_k}{l_k}}=(s+l_k)\cdots(s+1)/l_k !$.
The inequality is strict since there exists an
index $k$ such that $ m_k'>m_k$.
\end{proof}

\begin{Definition} \label{Dmpi}
For each choice of a permutation $\pi :\{1,\ldots,r\} \to
\{1,\ldots,r\}$, consider the degree vector $m^\pi$ defined
by
$$
m^\pi_{k} := \left(1+ \sum_{\pi(j)\geq \pi(k)}l_{j}\right)
d_{k} - l_{k}, \, \, k=1,\ldots, r.
$$
\end{Definition}

When all defects are zero, these are the vectors defined in \cite{StZe}
 yielding  determinantal Sylvester formulae and they also coincide with
those
defined in (\ref{syldetvec}) in Theorem \ref{Tdetmvec}.

\begin{lemma} \label{Lbase}
If $m\in\ZZ^r_{\geq 0}$ yields a Sylvester-type matrix, it is
possible to define a permutation $\pi:[r]\to[r]$ such that for
$i= \pi^{-1}(1)$ it holds that $m_i \ge m^\pi_{i}$.
Moreover $H^0 (\PP^{l_i}, m_i-pd_i )\neq 0$ where
$p \leq 1+\sum_{\pi(j) \in J} l_j$, for any subset
$J$ such that $\emptyset \neq J \subset \{2,\dots,r\}$.
\end{lemma}

\begin{proof} 
For $p=n+1,\nu=1$, a necessary condition is that $H^{n}(X, m-(n+1)d)=0$.
Hence, there exists $i \in [r] :
H^{l_i} (\PP^{l_i}, m_i-(n+1)d_i )=0 \Leftrightarrow
m_i-(n+1)d_i\ge -l_i \Leftrightarrow m_i\ge m^\pi_{i}$
by choosing  $\pi(i)=1$.
For any $p$ as in the statement,
$H^0 (\PP^{l_i}, m_i-pd_i ) \neq 0\Leftrightarrow m_i\ge pd_i$.
Since $m_i-(n+1)d_i\ge -l_i$ it suffices to prove
$(n+1)d_i -l_i\ge d_i (1+ \sum_{j\neq i} l_j ) \ge d_i p$.
The latter inequality is obvious for all $p$,
whereas the former
reduces to $l_i d_i\ge l_i$ which holds since $d_i \geq 1$.
\end{proof}

\begin{theorem}\label{sylvth}
A degree vector $m\in \ZZ^r_{\geq 0}$ gives a Sylvester-type matrix
iff there exists
a permutation $\pi$ such that $m_j \geq m^\pi_j$ for
$j=1,\ldots,r.$
Moreover, the smallest Sylvester matrix is attained among the vectors
$m^\pi$.
\end{theorem}
\begin{proof}
We prove the forward direction by induction on $k=1,\dots,r$.
Assume $m$ gives a Sylvester-type complex and
consider the necessary condition  $K_1(m)=H^0(X, m-d)$.
The base case $k=1$ was proven in Lemma \ref{Lbase}.
The inductive hypothesis for $k\in\{1,\dots,r-1\}$ specifies
which cohomologies vanish and which not, where
$m_{u}\ge m^{\pi}_{u},\, \pi(u) \le k$.
In particular, for all subsets $J$ such that
$\emptyset\neq J\subset \{1,\dots,r\}\setminus \{1,\dots,k\},\;
p=1+ \sum_{\pi(j) \in J} l_j,\;
p_0=p+l_v, $
for some $v$ such that $\pi(v) \le k$, we assume:
\begin{equation} \label{Estep}
H^{l_u} \left(\PP^{l_u}, m_u- p_0 d_u \right) =0,\;
H^{0} \left(\PP^{l_u}, m_u- p d_u \right) \neq 0.
\end{equation}
%
For the inductive step, we
exploit the necessary condition that $H^{p-1}(X, m-p d)=0$
for $p=1+ \sum_{\pi(j) > k} l_j$.
By (\ref{Estep}), there exists $i$  such that 
$H^{l_i} \left(\PP^{l_i}, m_i- p d_i \right) =0$. Then,
$m_i\ge pd_i-l_i=m^{\pi}_{i}$ where we define $\pi(i)=k+1$.
To complete the step, we show
$H^{0} \left(\PP^{l_u}, m_u- p d_u \right) \neq 0$
where $\pi(u) \le k+1,\,
p=1+ \sum_{\pi(j) \in J} l_j$, and any subset $J$ such that
$\emptyset\neq J\subset \{k+2,\dots,r\}$.
The non-vanishing of the cohomology is equivalent to
$m_u \ge p d_u$.
It suffices to prove
$m^\pi_{i} \ge d_i(1+ \sum_{\pi(j) > k+1} l_j)$.
By definition, this reduces to
$ -l_i + d_i l_{i} \ge 0\Leftrightarrow d_i\ge 1$.
The converse direction follows from analogous arguments as above.
The claim on minimality follows from Lemma \ref{Lstairs}.
\end{proof}

This gives an algorithm for finding the minimal Sylvester formulae
by testing at most $r!$ vectors $m^{\pi}$, which is implemented in
\Maple\  (Section~\ref{Scode}).
To actually obtain the square submatrix whose determinant is divisible by
the sparse resultant, it suffices to execute a rank test.
These matrices exhibit quasi-Toeplitz structure, implying
that asymptotic complexity is quasi-quadratic in the
matrix dimension \cite{EmPa02snap}.
Observe that $P_k(m^{\pi})\neq\emptyset$
because there exists $p\in\ZZ$  such that $p=1+\sum_{\pi(j)\le \pi(k)} l_j$
such that $m_k^{\pi} <d_k p= m_k^{\pi} +l_k$ for all $k$.

\begin{example} \label{Ex211_222}
Let $l=(2,1,1)$, $d=(2,2,2)$;
the degree of the resultant is $960$.
Let $\sigma=\pi^{-1}$ be the permutation inverse to $\pi$; then
the corresponding degree vector can be written as
$m^\pi_{\sigma(k)} := \left(1+ \sum_{j\geq k}l_{\sigma(j)}\right)
d_{\sigma(k)} - l_{\sigma(k)}$.
Here is a list of the $6= 3!$ degree vectors $m^{\pi}$,
among which we find the
smallest Sylvester matrix of row dimension $1080$, whereas
the sparse resultant's degree is 960.
Also shown are the
permutations $\sigma$ and the corresponding matrix dimensions.
The symmetry between the last two polynomials makes certain
dimensions appear twice.
$$
\begin{array}{rcrcc}
m^\pi=& (8,5,3) &\quad \sigma=&(1,2,3) & \quad 1080\times 1120 \\
& (8,3,5) & & (1,3,2)  & \quad 1080\times 1120\\
& (6,9,3) & & (2,1,3) & \quad 1120\times 1200\\
& (4,9,7) & & (2,3,1)  & \quad 1200\times 1440\\
& (6,3,9) & & (3,1,2)  & \quad 1120\times 1200 \\
& (4,7,9) & & (3,2,1) & \quad 1200\times 1440
\end{array}
$$
Our \Maple\  program, discussed in Section \ref{Scode},
enumerates 81 purely rectangular Sylvester matrices (none of
which is determinantal).
All Sylvester matrices not shown here have dimensions
$1260\times 1400$ or larger.
\end{example}
The map $K_1(m)\rightarrow K_0(m)$ is surjective, i.e.,
the matrix has at least as many columns as rows.
In searching for a minimal formula, we should reduce
$\dim K_0(m)$, i.e., the number of rows, since this defines
the degree of the extraneous factor in the determinant.
It is an open question whether $\dim K_0(m)$ reduces
iff $\dim K_1(m)$ reduces.
%
In certain system solving applications, the extraneous
factor simply leads to a superset of the common isolated
roots, so it poses no limitation.
Even if it vanishes identically, perturbation
techniques yield a nontrivial projection operator
\cite{DAnEmi01ams}.

It is possible to obtain a pure 
Sylvester-matrix
whose determinant equals the multihomogeneous
resultant when the complex has as only non-zero
terms
$K_1(m)=H^j(X, m-(j+1)d)^{\binom {n+1}{j+1}}$,
and $K_0(m)=H^j(X, m-jd)^{\binom {n+1} j}$
for any $j=0,\dots,n$.   But in this case we deduce  
from Theorem~\ref{twoterms} that
all defects vanish.
So (cf.\ \cite{StZe}) there exists a pure Sylvester-type
determinantal formula associated to a non negative degree vector
(i.e.\ $m_k \geq 0$, $k=1,\ldots,r$), or equivalently, for $j=0$.
Thus, if a determinantal data $(~l_1,\dots,l_r; d_1,\dots,d_r)$ admits
a degree vector yielding a two term Sylvester complex for some $j$, it
admits such a formula for $j=0$ as well.
Hence, concentrating on $j=0$ is not restrictive in the case
of determinantal complexes. 
The Sylvester-type matrices for positive $j$ correspond to
degree vectors with some negative entries, unlike the assumption
in \cite[p.118]{StZe}. We show such an example in the bilinear case
in \S \ref{Sexamples}. Thus,  the first part of conjecture \cite[Conj.3]{StZe}
can be true only for non-negative degree vectors.

\section{Pure B\'ezout-type formulae} \label{Sbezout}

In this section, we will study the following complexes:

\begin{Definition} \label{defbezfla}
A Weyman complex is of pure B\'ezout type if
$K_{-1}(m)=0$, $K_1(m)=H^{l_1+\cdots+l_r}\left(X, m-(n+1)d \right)$
and $K_0(m) =H^0(m)$.
\end{Definition}

Weyman complexes of pure B\'ezout type correspond
to generically surjective maps

\begin{equation} \label{EK10bezout}
H^{l_1+\cdots+l_r}\left(X, m-(n+1)d \right)\rightarrow
H^0(m) \rightarrow 0
\end{equation}
such that any maximal minor
is a nontrivial multiple of the multihomogeneous resultant.
In fact, we shall show that the only possible such formulae
are determinantal (i.e. $K_2(m)=0$).
We shall exhibit the corresponding differential
in terms of the Bezoutian and characterize the possible
degree vectors.
We show that there exists a pure B\'ezout-type formula
iff there exists a pure Sylvester formula. We remark that the
dimension of the matrix  with pure B\'ezout coefficients equals the
dimension
of the Sylvester matrix divided by $n+1$.
Now we can generalize results in \cite{ChtKap00}, \cite{Saxe97}
(cf.\ Section \ref{Sintro}).
\begin{theorem}
There exists a determinantal formula of pure
B\'ezout type iff for all $k$ either $l_k =1$ or $d_k =1$,
i.e.\  all defects vanish.
\end{theorem}
\begin{proof}
This is just a special instance of Theorem~\ref{twoterms}, where
the only non-zero cohomologies in the complex correspond to
$p_1=0, p_2 = n+1.$
\end{proof}

Let us study degree vectors yielding pure B\'ezout formulae, which will
then provide the smallest determinantal formulae in case all defects
vanish.

\begin{Definition} \label{Dmbezout}
For each choice of a permutation $\pi :\{1,\ldots,r\} \to
\{1,\ldots,r\}$, let us define a degree vector
$$
m^\pi_{k} := - l_{k} +
d_{k} \sum_{\pi(j)\geq \pi(k)}l_{j},\quad \, k=1,\dots, r.
$$
\end{Definition}

When all $\delta_k=0$ these are precisely the vectors defined
in (\ref{bezdetvec}) in Theorem \ref{Tdetmvec2}. Note that
our assumptions in
(\ref{EK10bezout}) imply
$H^{l_j}\left(\PP^{l_j}, m_j-(n+1)d_j \right)\neq 0$,
$H^0\left(\PP^{l_j}, m_j\right)\neq 0$.
\begin{lemma} \label{LboundmB}
The existence of any pure B\'ezout formula implies
$0\le m_j<(n+1)d_j-l_j$, for all $j$.
\end{lemma}
In fact, the $m^\pi$ of Def.~\ref{Dmbezout}
satisfy these constraints for all permutations $\pi$.

\begin{proof}
The non negativity of all $m_j$ is deduced from the fact
that $H^0(X, m) \not= 0.$ On the other side, the non vanishing of
$H^n(X, m-(n+1)d)$  implies the other inequality.
\end{proof}

\begin{lemma} \label{LbaseBez}
If $m\in\ZZ^r$ yields a pure B\'ezout-type complex,
then there exists a permutation $\pi:[r] \rightarrow [r]$
such that  $i = \pi^{-1}(1)$ verifies $\, m_i \ge m^\pi_{i}$ and
\begin{eqnarray*}
H^{l_i}\left(\PP^{l_i}, m_i-(q+l_i+\nu)d_i\right)=0,\quad \nu=0,-1,\\
H^{0}\left(\PP^{l_i}, m_i -q d_i\right)\neq 0,
\quad q=\sum_{j\in J} l_j,
\end{eqnarray*}
for any $J\subset\{1,\dots,r\}\setminus\{i\}, J \neq \emptyset$.
\end{lemma}

\begin{proof}
Since $H^n(X, m-nd) =0$, there exists an index
$ i\in\{1,\dots,r\}$ such that $m_i - n d_i \geq -l_i$. It is enough
to define $\pi(i) = 1.$
\end{proof}

\begin{theorem} \label{LinduceBez}
If $m\in\ZZ^r$ yields a pure B\'ezout-type complex,
it is possible to find a permutation $\pi$ such that
the degree vector $m$ verifies
$m_i \ge m^{\pi}_i$ for all $i=1,\dots,r$.
\end{theorem}

\begin{proof} 
We use induction;
the base case follows from Lemma \ref{LbaseBez}.
The inductive hypothesis, for $k\in\{1,\dots,r-1\}$, is:
there exists a subset $ U\subset\{1,\dots,r\},\, |U|=k$,
such that $\pi(u) \le k$, $m_u \ge m^\pi_{u}$ for all $u\in U$ and
\begin{eqnarray}\label{EhypoBez}
H^{l_u}\left(\PP^{l_u}, m_u-(q+l_u+\nu)d_u\right)=0,\quad \nu=0,-1,\\
H^{0}\left(\PP^{l_u}, m_u -q d_u\right)\neq 0,
 \quad q=\sum_{j\in J} l_j, \nonumber
\end{eqnarray}
for all $ J\subset\{1,\dots,r\}\setminus U, J\neq \emptyset$.
Now the inductive step:
The hypothesis on $K_0$ implies
$H^{p}(X, m-pd)=0$ for $p=\sum_{j\not\in U} l_j$.
Considering the inequality in (\ref{EhypoBez}) for $q=p$,
there exists $ i\in [r]\setminus U$ such that 
$H^{l_i}\left(\PP^{l_i}, m_i- p d_i\right)=0\Leftrightarrow
m_i+l_i\ge p d_i$ i.e.\  $m_i \ge m^\pi_i$ for
$\pi(i)=k + 1$
because $j\not\in U\Leftrightarrow \pi(j)\ge k + 1$.
It suffices now to extend (\ref{EhypoBez})
for $q'=\sum_{j\in J'} l_j$ where
$ \emptyset\neq J'\subset [r] \setminus (U\cup\{i\})$.
First, $m_i+l_i\ge p d_i \ge (q'+l_i) d_i$ implies
the equations below.
Second, $m_i\ge -l_i+pd_i=(p-l_i)d_i + l_i(d_i-1)\ge q'd_i$
yields the inequality, so
\begin{eqnarray}\label{EsteBez}
H^{l_i}\left(\PP^{l_i},m_i-(q'+l_i+\nu)d_i\right)=0,\; \nu=0,-1.\\
H^{0}\left(\PP^{l_i},m_i -q' d_i\right)\neq 0.\nonumber
\end{eqnarray}
Then, $m$ satisfies the hypothesis $K_{-1}(m)=0$
for $\nu=-1$ because every summand in $K_{-1}(m)$ contains some
cohomology as in (\ref{EsteBez}).
Since $p\ge 0\Rightarrow q=p-\nu\ge 1$ no
summand has only zero cohomologies.
By Lemma \ref{LboundmB} and (\ref{EsteBez}) for $\nu=0$,
$m$ gives
$K_0(m)=H^0(X, m)$ because $H^{l_i}(\PP^{l_i},m_i-pd_i)=0$
for any $p,i$.
\end{proof}

\begin{theorem} \label{Tbezout}
A pure B\'ezout and generically surjective formula exists for some
vector $m$ iff it equals $m^\pi$ of Definition
\ref{Dmbezout},
for some permutation $\pi$, and all defects are zero.
\end{theorem}

\begin{proof} 
It suffices to consider $K_1(m) =H^n(X, m-(n+1)d)$;
it is nonzero by Lemma \ref{LboundmB}.
We prove by induction that $m=m^\pi$ by using the fact
that all other summands
in (\ref{EKnu_coH}) for $K_1(m)$ vanish.
For $H^0(X, m-d)$ to vanish, there must exist
$i\in [r]$ such that $H^0(\PP^{l_i}, m_i-d_i)=0\Leftrightarrow
m_i < d_i$.
Hence we need to define $\pi(i)=r$ because
$\pi(i)< r \Rightarrow m^\pi_i \ge -l_i+d_i (l_i+1)=d_i + l_i(d_i-1)
\ge d_i$.
Moreover, $m^\pi_i=-l_i+d_il_i<d_i \Leftrightarrow \delta_i=0$
by Lemma \ref{Ldefect}.

There is a unique integer in $[m_i^\pi,d_i)$ because
$m_i^\pi+1\ge d_i\Leftrightarrow -l_i+d_il_i+1\ge d_i
\Leftrightarrow (l_i-1)(d_i-1)\ge 0$.
Hence $m_i=d_i-1<d_i(q+1)$ for any $q\ge 0$, therefore
$H^0(\PP^{l_i}, m_i-(q+1)d_i)=0$.
Furthermore, for $q\ge l_i$,
$H^{l_i}(\PP^{l_i}, m_i-(1+q)d_i) \neq 0\Leftrightarrow m_i+l_i<(q+1)d_i
\Leftrightarrow l_i-1<qd_i$ which holds.
This proves the inductive basis.
The inductive hypothesis is: for all $ u\in U\subset [r]$,
where $|U|=k,\, \pi(u) > r-k$, then $\delta_u=0$,
$m_u=m^\pi_u$ and
\begin{equation}\label{EhypoK1bez}
\! H^0(\PP^{l_u}, m_u- (1+\sum_{\pi(j)<\pi(u)} l_j )d_u)=0\neq
H^{l_u}(\PP^{l_u}, m_u- (1+ \sum_{j\in J} l_j )d_u),
\end{equation}
for all $J$ such that $U\subset J\subset [r]$.
For the inductive step, consider that $H^q(X, m-(1+q)d)$ must vanish
for $q=\sum_{j\in U} l_j$.
None of its summand cohomologies
$H^{l_u}\left(\PP^{l_u}, m_u- (1+ q)d_u\right)$ vanish due to
the last inequality.
So there exists  $i$ such that $ H^0(\PP^{l_i}, m_i - (1+q)d_i)=0\Leftrightarrow
m_i<(1+q)d_i$.

Hence $\pi(i)=r-k$ so that
$m_i = m^\pi_i = -l_i + d_i\sum_{\pi(j)\ge r-k} l_j < (1+q)d_i
\Leftrightarrow -l_i + d_il_i < d_i \Leftrightarrow \delta_i =0$
by Lemma \ref{Ldefect}.
No larger $m_i$ works because $m_i^\pi$ is the maximum integer
strictly smaller than $(1+q)d_i$.
And $\pi(i)<r-k$ would make $m_i$ too large.
Now extend the inequality
(\ref{EhypoK1bez}) to $J'$ where $(U\cup\{i\})\subset J'$
and observe
$m_i^\pi < d_i\sum_{\pi(j)\ge r-k} l_j < d_i (1+\sum_{j\in J'}l_j)$.

The hypothesis is proven for all $U\subset [r]$,
including the case $|U|=r$.
For the converse, assume there exists a permutation $\pi$
such that $ m=m^\pi$ and all defects vanish.
Then $K_0(m), K_{1}(m)$ satisfy all
conditions for a pure B\'ezout formula.
Furthermore, $K_{-1}(m)=0$, hence the
formula is generically surjective.
\end{proof}

The condition $K_2(m)=0$, which yields a square matrix, is obtained by
the hypothesis of a pure B\'ezout and generically surjective formula;
i.e., there is no rectangular surjective pure B\'ezout formula.

\begin{corollary} \label{Cbezout}
If a generically surjective formula is of pure B\'ezout type,
then it is determinantal.
Furthermore, for any permutation $\pi$, the matrix is of the
same dimension, i.e.\ $\dim K_0(m) = \deg R / (n+1).$
\end{corollary}

\subsection{Explicit B\'ezout-type formulae} \label{SexplBez}

We start by defining a dual permutation $\pi'$ to any permutation $\pi$.

\begin{Definition} \label{Dcdg_pprim}
For any permutation $\pi : [r]\rightarrow [r]$, define a new
permutation $\pi' : [r]\rightarrow [r]$
by $\pi'(i)=r+1-\pi(i)$.
\end{Definition}

\begin{lemma} \label{Lmsumrho}
Assuming all defects are zero,
$m^\pi+ m^{\pi'}=\rho$ for any permutation $\pi : [r]\rightarrow [r]$,
where $\rho\in\NN^r$ is the critical vector of Definition \ref{Dcrit}.
\end{lemma}

\begin{proof} 
$m^\pi_i+m^{\pi'}_i=d_i (n+l_i)-2l_i$ for all $i$
because the sum in the parenthesis includes
$\{l_j:\pi(j)\ge\pi(i)\} \cup \{l_j:\pi'(j)\ge\pi'(i)\}$,
and latter set is $\{l_j:\pi(j)\le\pi(i)\}$.
So $m^\pi_i+m^{\pi'}_i = d_il_i-l_i+\rho_i-d_i+1
= \rho_i+d_i(l_i-1)-(l_i-1)=\rho_i+(l_i-1)(d_i-1)=\rho_i$
because of the zero defects.
\end{proof}

Denote by $x_i$ (resp.\ $x_{ij}$) the $i$-th variable group
(resp.\ the $j$-th variable in the group), $i\in [r],
j=0,\dots,l_i$.
Introduce $r$ new groups of variables $y_{i}$ with the same
cardinalities and denote by $y_{ij}$ their variables.

Given a permutation $\pi$, let the associated {\em Bezoutian}
be the polynomial $B^\pi(x,y)$ obtained as follows:
First dehomogenize the polynomials by setting $x_{i0}=1,\,
i=1,\dots,r$;
the obtained polynomials are denoted by $f_0,\ldots, f_n$.
Second, construct the $(n+1)\times(n+1)$ matrix with
$j$-th column corresponding to polynomial $f_j,\, j=0,\dots,n,$
and whose $x_{ij}$ variables are gradually substituted, in successive
rows, by each respective $y_{ij}$ variable.
This construction is named after B\'ezout or Dixon and is
well-known in the literature, e.g.\ \cite{CaMo95}, \cite{EmiMou99jsc}.
A general entry is of the form
\begin{eqnarray} \label{EentryBez}
\lefteqn{f_j (y_{\sigma(1)},\dots,y_{\sigma(k-1)},y_{\sigma(k)1},
 \dots,y_{\sigma(k)t},}\\
& x_{\sigma(k)(t+1)},\dots,x_{\sigma(k)l_{\sigma(k)}},x_{\sigma(k+1)},
 \dots,x_{\sigma(r)} )\nonumber
\end{eqnarray}
where $\sigma:=\pi^{-1}, k=0,\dots,r, t=1,\dots,l_k$.
There is a single first row for $k=0$, containing all the polynomials
in the $x_{ij}$ variables, whereas the last row has the same
polynomials
with all variables substituted by the $y_{ij}$.
All intermediate rows contain the polynomials in a subset of the
$x_{ij}$
variables, the rest having been substituted by each corresponding
$y_{ij}$.
The number of rows is $1+\sum_{j \in [r]} l_j=1+n$.
Lastly, in order to obtain $B^\pi(x,y)$,
we divide the matrix determinant by
\begin{equation}\label{EprodDif}
\prod_{i =1}^r \prod_{j=1}^{l_i} (x_{ij}-y_{ij}).
\end{equation}

\begin{example} \label{Ex12_21}
Let $l=(1,2), d=(2,1)$.
If $\pi=(12), \pi'=(21)$, then $m^\pi=(5,0), m^{\pi'}=(1,1)$.
For both degree vectors, the matrix dimension is 6. 
To obtain $B^\pi(x,y)$ we construct a $4\times 4$ matrix
whose $j$-th column contains
$f_j(x_{1},x_{21},x_{22}), f_j(y_{1},x_{21},x_{22})$,
$f_j(y_{1},y_{21},x_{22})$, $ f_j(y_{1},y_{21},y_{22})$,
for $j=0,\dots,3$.
Here $x_1$ (and $y_1$) is a shorthand for $x_{11}$ (and $y_{11}$).
Then $B^\pi(x,y)$ contains the following monomials
in the $x_i$ and $y_i$
variables respectively, 6 in each set of variables:
$1, x_1, x_{21}, x_{22}, x_1x_{21}, x_1x_{22},$
$1, y_1, y_1^{2}, y_1^{3}, y_1^{4}, y_1^{5}.$
So the final matrix is indeed square of dimension 6.
More details on this example are provided in Section \ref{Scode}.
\end{example}
\begin{lemma}  \label{Lboundab}
Let $B^\pi(x,y) = \sum b_{\alpha \beta} x^\alpha y^\beta$
where $\alpha = (\alpha_{ij}), \beta = (\beta_{ij}) \in\ZZ^n$,
$i= 1, \dots, r, j= 1, \dots, l_i$.
Set
$\alpha_i=$ $\sum_{j=1,\dots,l_i}\alpha_{ij}$,
$\beta_i= \sum_{j=1,\dots,l_i}\beta_{ij}$,
for all $ \alpha, \beta$.
Then,
$0\le \alpha_i\le m^{\pi'}_i, 0\le \beta_i \le m^\pi_i$
and $0\le \alpha_i+\beta_i \le \rho_i,\, i=1,\dots,r$.
\end{lemma}

\begin{proof} 
By Lemma \ref{Lmsumrho} it suffices to bound $\alpha_i,\beta_i$.
But $\alpha_i$ is the degree of the $x_i$ in the determinant
decreased by $l_i$ in order to account for the division by
(\ref{EprodDif}).
The former equals the product of $d_i$ with the number of rows
where an $x_{ij}$ variable appears for any $j\in [1,l_i]$.
These are the first row, the rows where $y_j$ are introduced for
$j\in\{\sigma(1), \dots, \sigma(k-1)\}$ such that $\sigma(k)=i$,
and another $l_i-1$ rows when $\sigma(k) = i$.
The condition on $j : \pi(j)<\pi(i)$ is equivalent to
$r+1-\pi(j)>r+1-\pi(i)$, hence
$
\alpha_i\le -l_i+d_i \sum_{\pi'(j) \ge \pi'(i)} l_j =m_i^{\pi'}.
$
Similarly, we prove the upper bound on $\beta_i$.
The rows containing $y_{ij}$ for some $j\in [1,l_i]$ are those
where $j\in\{\sigma(k+1), \dots, \sigma(r)\} : \sigma(k)=i$,
another $l_i-1$ rows when $\sigma(k) = i$, and the last row.
Now, $\pi(j)\ge k+1>k=\pi(i)$, so
$
\beta_i\le -l_i+d_i \sum_{\pi(j) \ge \pi(i)} l_j =m_i^{\pi}.
$
Clearly $\alpha_i, \beta_i \ge 0$.
\end{proof}
For generic polynomials, the upper bounds of $\alpha_i, \beta_i$ are
attained.
The lemma thus gives tight bounds on the support of the Bezoutian.
\begin{theorem}\label{explicitbez}
Assume all defects are zero and $B^\pi(x,y)$ is defined as above.
For any $\pi$, $(b_{\alpha \beta})$ is a square
matrix of dimension
$$
\dim K_0(m) = {\binom l{l_1,\dots,l_r}} d_1^{l_1}\cdots d_r^{l_r}
= \frac{ \deg R }{ (n+1) }.
$$
Furthermore, $\det (b_{\alpha \beta})=R(f_0,\dots,f_n)$.
\end{theorem}
\begin{proof}
First, we show that $(b_{\alpha \beta})$ is square of the desired size.
The dimensions are given by the number of exponent vectors
$\alpha,\beta$ bounded by Lemma \ref{Lboundab} which
are exactly $\dim K_0(m^{\pi'}), \dim K_0(m^{\pi})$ respectively.
Both $m^\pi, m^{\pi'}$ are determinantal, hence
both of these numbers are equal to $\deg R / (n+1)$,
by Theorem \ref{Tbezout} and Corollary \ref{Cbezout}.
$R(f_0,\ldots,f_n)$ divides
every nonzero maximal minor of the matrix  $(b_{\alpha \beta})$;
cf.\ \cite{CaMo95}, \cite[Thm.~3.13]{EmiMou99jsc}.
Since any nonzero proper minor has degree
$< \deg R$,  the determinant of the
matrix  $(b_{\alpha \beta})$ is nonzero and equals the resultant.
\end{proof}

Note that there is not a unique choice of higher
differentials in the Weyman complexes.
We could chase the arrows in a resultant spectral sequence as in
\cite[Ch.~2, Prop.~5.4]{GKZ}
to show that the matrix we propose comes from 
the explicitization of one possible choice.
We have followed instead the more direct route based on the above property
of the Bezoutian in \cite{CaMo95}, which uses more elementary tools.

\section{Two examples}   \label{Sexamples}

\subsection{The bilinear system}   \label{Sbilinear}

The generic system of 3 bilinear polynomials is
$$
\begin{array}{c}
f_0 = a_0 + a_1 x_1+ a_2x_2+ a_3 x_1x_2,\\
f_1= b_0 + b_1 x_1+ b_2x_2+ b_3 x_1x_2,\\
f_2= c_0 + c_1 x_1+ c_2x_2+ c_3 x_1x_2,\\
\end{array}
$$
and has type $(1,1;1,1)$.
The degree of the resultant in the coefficients is $3{\binom 2{1,1}} = {6}$.
We shall enumerate all 14 possible determinantal formulae in the
order of decreasing matrix dimension, from 6 to 2, and shall
make the corresponding maps explicit. This study goes back to the
pioneering work of Dixon \cite{dix}.

For $\pi=(1,2)$, Definition \ref{Dmpi} yields $m=(2,1)$ and the complex
is $0\to K_1=H^0(1,0)^{{\binom 31}}\to K_0=H^0(2,1)\to 0$.
The corresponding determinantal pure Sylvester matrix is,
when transposed, equal to
$$
\left[ \begin{array}{cccccc}
a_0 & a_1 & a_2 & a_3 & 0   & 0\\
b_0 & b_1 & b_2 & b_3 & 0   & 0\\
c_0 & c_1 & c_2 & c_3 & 0   & 0\\
0   & a_0 & 0   & a_2 & a_1 & a_3\\
0   & b_0 & 0   & b_2 & b_1 & b_3\\
0   & c_0 & 0   & c_2 & c_1 & c_3\\
\end{array}\right],
$$
with rows corresponding to the input polynomials and the same
set multiplied by $x_1$, whereas the columns are indexed by
$1, x_1, x_2, x_1x_2, x_1^2, x_1^2x_2$.
By symmetry, another formula is possible by interchanging the
roles of $x_1,x_2$.
Further formulae are obtained by taking the transpose of these
2 matrices, namely with $m=(-1,0)$, where the complex is
$H^2(-4,-3)=H^0(2,1)^*\to (H^2(-3,-2))^3=(H^0(1,0)^*)^3$,
and with $m=(0,-1)$.
Recall the definition of duality from equation (\ref{Serre}).
Sylvester maps, as well as other types of maps, are further
illustrated in Section \ref{Sexample}.

There are additional determinantal
Sylvester formulae corresponding to $m=(2,-1)$ and $m=(-1,2)$.
Their matrices contain the $f_i$ and the $f_i$ multiplied by
$x_1^{-1}$ or by $x_2^{-1}$.
In the former case, the complex is $H^1(0,-3)^3 =
(H^0(0)\otimes H^0(1)^*)^3 \to H^1(1,-2)^3 = (H^0(1)\otimes H^0(0)^*)^3$,
the transposed matrix is
$$
\left[ \begin{array}{cccccc}
a_0 & a_1 & a_2 & a_3 & 0   & 0\\
b_0 & b_1 & b_2 & b_3 & 0   & 0\\
c_0 & c_1 & c_2 & c_3 & 0   & 0\\
a_1 & 0   & a_3 & 0   & a_0 & a_2\\
b_1 & 0   & b_3 & 0   & b_0 & b_2\\
c_1 & 0   & c_3 & 0   & c_0 & c_2\\
\end{array}\right],
$$
and the columns are indexed by
$1, x_1, x_2, x_1x_2, x_1^{-1}, x_1^{-1}x_2$.
This construction can be verified by hand calculations.

For $m=(1,1)$ the complex becomes
$0\to K_1=H^0(0,0)^{{\binom 31}}\oplus H^2(-2,-2)\to K_0=H^0(1,1)\to 0$.
The matrix is square, of dimension equal to 4, and hybrid.
We compute the two maps by hand; for a larger example see 
Section~\ref{Sexamples}.
The foundations for constructing such matrices can be found in \cite{WeZe}.
The transposed $4\times 4$ determinantal formula is written
as follows, by using brackets:
$$
\left[ \begin{array}{cccc}
a_0 & a_1 & a_2 & a_3 \\
b_0 & b_1 & b_2 & b_3 \\
c_0 & c_1 & c_2 & c_3 \\
{}[012] & [013] & [032] & -[123]
\end{array}\right],
 \; \mbox{ where } \quad [ijk]= \det \left[ \begin{array}{ccc}
a_i&a_j&a_k\\ b_i&b_j&b_k\\ c_i&c_j&c_k\end{array}\right].
$$
The matrix rows contain the $f_i$ and a rational  multiple  of the affine toric Jacobian,
whereas the columns are indexed by $1,x_1,x_2,x_1x_2$.
This formula is obtained in \cite{CaDiSt98} in a more general toric setting.
An analogous $4\times 4$ matrix corresponds to $m=(0,0)$.

There are 4 ``partial B\'ezout'' determinantal formulae of dimension
$3\times 3$ for $m=(-1,1), (1,-1)$ and for $m=(2,0), (0,2)$. We
omit the details of the computation.
In the first case, the complex is
$H^2(-4,-2)=H^0(2,0)^* \to H^1(-2,0)^3= (H^0(0)^*\otimes H^0(0))^3$,
and a choice of the matrix is, in terms of brackets,
$$
\left[ \begin {array}{ccc}
[-02] & [-03]+[-12] & [-13]\\
{}[0\!-\!2] & [0\!-\!3]+[1\!-\!2] & [1\!-\!3]\\
{}[02-] & [12-]+[03-] & [13-]
\end{array}\right],
 \; \mbox{ where } \quad [ij-]= \det \left[ \begin{array}{cc}
a_i&a_j\\ b_i&b_j\end{array}\right],
$$
and analogously for the $2\times 2$ brackets $[i-k], [-jk]$.
The columns of this resultant matrix are indexed by $1,x_1,x_1^2$,
which is the support of the 3 Bezoutian polynomials filling in the rows.
In particular, these polynomials are defined for $\{i,j,k\}=\{0,1,2\}$
in the standard way:
$$
B_k= \det\left[ \begin {array}{cc}
f_i(x_1,x_2) & f_i(x_1,y_2) \\
f_j(x_1,x_2) & f_j(x_1,y_2) \\
\end{array}\right] \; / \; (x_2-y_2).
$$

For $m=(1,0)$ the complex becomes
$0\to K_1=H^2(-2,-3)^{{\binom 33}}=H^0(0,1)^*\to K_0=H^0(1,0)\to 0$.
The corresponding determinantal pure B\'ezout-type formula
is obtained from the Bezoutian polynomial
$$
B= \det\left[ \begin {array}{ccc}
f_0(x_1,x_2) & f_0(y_1,x_2) & f_0(y_1,y_2)\\
f_1(x_1,x_2) & f_1(y_1,x_2) & f_1(y_1,y_2)\\
f_2(x_1,x_2) & f_2(y_1,x_2) & f_2(y_1,y_2)\\
\end {array}\right] \; / \; (x_1-y_1)(x_2-y_2),
$$
supported by $\{1,x_2\},\{1,y_1\}$.
The resultant matrix is given in terms of brackets as follows:
$$
\left[ \begin {array}{cc}
 [123] & [023] \\\noalign{\bigskip}
-[103] & [012]
\end{array}\right].
$$

\subsection{A hybrid determinantal formula}   \label{Sexample}

Assume $l=(3,2), d=(2,3)$. We present explicit
formulae which can be extrapolated in general, giving
an answer to the problem stated in \cite[p.~578]{WeZe}. We plan to
carry this extensively in a future work, but we include here the
example without proofs as a hint for the interested reader.
Our \Maple\
program enumerates 30 determinantal vectors $m$, among which
we find $m^{(2,1)}=(3,13)$ according to Theorem~\ref{Tdetmvec}.

The minimal matrix dimension is $1320$ and
is achieved at $m=(6,3)$ and $(2,12)$.
In both cases, $P_2 (m)= \emptyset$, whereas $P_1(6,3)=\{4\}$
and $P_1(2,12)=\{2\}$.
This shows that the minimum matrix dimension may occur for some
empty $P_k$, contrary to what one may think.

Moreover, the degree of the sparse resultant is
$6{\binom 5{3,2}}2^3 3^2=4320$.
Since $1320$ does not divide $4320$, the minimal matrix is not of
pure B\'ezout type; it is not of pure Sylvester type either.
To specify the cohomologies and the linear maps that make the
matrix formula explicit we compute, for the degree vector $m=(6,3)$
and $p=1,\dots,6$ the different values of $m-pd$:
$(4,0)$, $(2,-3)$, $(0,-6)$, $(-2,-9)$, $(-4,-12)$, $(-6,-15)$.
The complex becomes
$K_2=0\rightarrow K_1\rightarrow K_0\rightarrow K_{-1}=0$,
with nonzero part
\begin{eqnarray*}
\lefteqn{H^0(4,0)^{{\binom 61}}\oplus H^2(0,-6)^{\binom63}
\oplus H^5(-6,-15)^{{\binom 66}} \rightarrow} \\
&\rightarrow H^0(6,3)^{{\binom 60}}
\oplus H^2(2,-3)^{{\binom 62}}\oplus H^5(-4,-12)^{{\binom 65}},
\end{eqnarray*}
where we omitted the reference to the space $X= \PP^{3}
\times \PP^{2}$ in the notation of the cohomologies.
Then $\dim K_1= 210+200+910=1320=840+150+330=\dim K_0$.
By a slight abuse of notation, let $\delta_{\alpha,\beta}$ stand for
the restriction of the above map to $H^{\alpha}\rightarrow H^{\beta}$.
Then $\delta_{02}=\delta_{05}=\delta_{25}=0$ by \cite[Prop.~2.5]{WeZe}
and it suffices to study the maps below, of which the first 3
are of pure Sylvester type by \cite[Prop.~2.6]{WeZe} and the last 3
are of pure B\'ezout type as those of Section \ref{Sbezout}.
These maps can be simplified using the dual cohomologies:
$$
H^{j}(\PP^{l_k}, m_k-pd_k)=H^{l_k-j}(\PP^{l_k}, (\rho_k-m_k)-(n+1-p)d_k)^*,
$$
where $\rho$ is the critical vector of Definition \ref{Dcrit}. So, we
have maps
\begin{eqnarray*}
\delta_{00}:& H^0(4,0)^{6} \rightarrow H^0(6,3)\\
\delta_{22}:& (H^0(0)\otimes H^0(3)^*)^{{\binom 63}}
 \rightarrow (H^0(2)\otimes H^0(0)^*)^{{\binom 62}} \\
\delta_{55}:& H^0(2,12)^* \rightarrow (H^0(0,9)^*)^{6}\\
\delta_{20}:& (H^0(0)\otimes H^0(3)^*)^{{\binom 63}}
 \rightarrow H^0(6,3)\\
\delta_{50}:& H^0(2,12)^* \rightarrow H^0(6,3) \\
\delta_{52}:& H^0(2,12)^* \rightarrow (H^0(2) \otimes H^0(0)^*)^{{\binom
62}}\\
\end{eqnarray*}
The resultant matrix (of the previous map in the natural
monomial bases)
has the following aspect, indicated
by the row and column dimensions:
$$
\begin{array}{lcr}
 & \, 840 \quad 150 \quad 330 & \\
\begin{array}{l} 210\\ 220\\ 910\\ \end{array}
&
\left[ \begin{array}{ccc}
\delta_{00}\quad & 0		&\quad 0\\
\delta_{20}\quad & \delta_{22}	&\quad 0\\
\delta_{50}\quad & \delta_{52}	&\quad \delta_{55}
\end{array} \right]
&
= \left[ \begin{array}{ccc}
S_{00}          & 0             & 0\\
B^{x_2}_{20}    & \delta_{22}   & 0\\
B_{50}          & B^{x_1}_{52}  & S_{55}^T
\end{array} \right]
\end{array}
$$
where $S_{ij}, B_{ij}, B^{x_k}_{ij}$ stand for pure Sylvester and
B\'ezout blocks,
the latter coming from a Bezoutian with respect to variables $x_k$ for
$k=1,2$,
and
$S_{55}^T$ represents a transposed Sylvester matrix, corresponding to
the dual
of the Sylvester map $H^0(0,9)^6\rightarrow H^0(2,12)$.

Let us take a closer look at $\delta_{22}$, which denotes both
the map and the corresponding matrix.
Let $\alpha\in\NN^{l_1}, \beta\in\NN^{l_2}$, be the degree vectors of
the
elements of $H^0(2), H^0(3)^*$ respectively, thus $|\alpha|\le 2,
|\beta|\le 3$.
Let $I,J\subset \{0,\dots, 5\}, |I|=3, |J|=2$ express the chosen
polynomials
according to the cohomology exponents.
Then the entries are given by
$$
\delta_{22}(x_1^{\alpha}\otimes T_J, 1\otimes S_I^{\beta}) =
\left\{ \begin{array}{l} 0, \mbox{ if } J\not\subset I,\\
 \mbox{coef}(f_k) \mbox{ of } x_1^{\alpha}x_2^{\beta},
 \mbox{ if }I\setminus J=\{k\},
\end{array} \right.
$$
where $T_J\in H^0(0)^*, S_I^\beta$ are elements of the respective dual
bases of monomials.
We expect such a construction to be generalizable, but
such a proof would be part of future work.

Now take the B\'ezout maps:
The matrix entries are given in (\ref{EentryBez})
for $\sigma=$ $(2,1)$:
the entry $(i,j)$, $i,j\in\{0,\dots,5\}$ contains
$f_j(x^{(1)}$, $\dots$, $x^{(5-i)}$, $y^{(6-i)}$, $\dots,y^{(5)}),$
where each $x^{(i)}$ is a leading subsequence of
$x_{11}$, $x_{12}$, $x_{13}$, $x_{21}$, $x_{22}$;
similarly with the new variables $y^{(i)}$.
The degree of the determinant, i.e.\ the Bezoutian,
is $6,3,2,12$ in $x_1,x_2,y_1,y_2$ respectively and these
coefficients fill in  the matrix $B_{50}$.
For the B\'ezout block $B_{52}^{x_1}$, consider ``partial"
Bezoutians defined from the 6 polynomials with the exception of
those indexed in $J$, where $J, I$ are as above.
Only the $x_1$ variables are substituted by new ones, thus yielding
a $4\times 4$ matrix.
For $B_{20}^{x_2}$, take all polynomials indexed in $I$ and
develop the Bezoutian with new variables $y_2$ from a
$3\times 3$ matrix.
Hence the entries of the B\'ezout blocks have,
respectively, degree $6,4,3$ in the coefficients of the $f_i$.

\section{Implementation} \label{Scode}

\begin{table}\begin{center}\begin{tabular}{|l|p{11truecm}|}
\hline
routine  & function\\
\hline \hline
{\tt comp\_m} & compute the degree vector $m$ by some specified formula \\
\hline
{\tt allDetVecs} & enumerate all determinantal formulae \\
\hline
{\tt allsums} & compute all possible sums of the $l_i$'s adding
   to $q\in\{0,\dots,\sum_{i=1}^r l_i\}$.\\
\hline
{\tt coHzero} & test whether $H^q(X, m-pd)$ vanishes\\
\hline
{\tt coHdim} & compute the dimension of $H^q(X, 
m-pd)$\\
\hline
{\tt dimKv} & compute the dimension of $K_{\nu}$ i.e.\  of the
corresponding matrix\\
\hline
{\tt findBez} & find all $m$-vectors yielding a pure B\'ezout-type
formula \\
\hline
{\tt findSyl} & find all $m$-vectors yielding a pure Sylvester-type
formula\\
\hline
{\tt minSyl} & find all $m^{\pi}$-vectors yielding a pure Sylvester-type
formula\\
\hline
{\tt hasdeterm} & test whether a determinantal formula exists\\
\hline
\end{tabular}
\caption{The main functionalities of our software\label{Tmaple}.}
\end{center} \end{table}

We have implemented on {\sc Maple~V} routines for the
above operations, including those in Table \ref{Tmaple}.
They are illustrated below and are
available in file {\tt mhomo.mpl} through:\\
{\tt http://www.di.uoa.gr/\~{}emiris/index-eng.html}.

\begin{oneshot}{Example \protect\ref{Ex12_23} (Cont'd)}
Recall that $l=(1,2), d=(2,3)$ and let $m=(6,3)$:\\
{\tt
> Ns:=vector([1,2]): Ds:=vector([2,3]):\\
> summs:=allsums(Ns):\\
> hasdeterm(Ns,Ds,vector([6,3]),summs);\\
}
\centerline{\em true}
{\tt > dimKv(Ns,Ds,vector([6,3]),summs,1);\\ }
\centerline{88}
{\tt > dimKv(Ns,Ds,vector([6,3]),summs,0);\\ }
\centerline{88}
\end{oneshot}

\bigskip

\begin{oneshot}{Example \protect\ref{Ex211_222} (Cont'd)}
Recall that $l=(2,1,1), d=(2,2,2)$, then $\delta=(1,0,0)$.
The \Maple\  session first computes all 81 pure Sylvester formula
by searching the appropriate range of 246 vectors.
The smallest formulae are shown.\\
{\tt
> Ns:=vector([2,1,1]):Ds:=vector([2,2,2]):\\
> minSyl(Ns,Ds):
}
\begin{center}{\em list\ of\ minimal\ S-matrices:\ m-vector\ and\ K1,K0-dims}\end{center}
\begin{eqnarray*}
& [[8,\,5,\,3,\,1120,\,1080],\,[8,\,3,\,5,\,1120,\,1080],\,[6,\,9,\,3,\,1200,\,1120],\\
& [6,\,3,\,9,\,1200,\,1120],\,[4,\,9,\,7,\,1440,\,1200],\,[4, \,7, \,9, \,1440, \,1200]] 
\end{eqnarray*}
{\tt
> allSyl:=findSyl(Ns,Ds):\\
{\small
Search of degree vecs from [4,3,3] to [8,9,9].\\
First array [4,7,9]: dimK1=1440, dimK0=1200,\\
dimK(-1)=0(should be 0).\\
}
}
\centerline{{\em\#pure-Sylvester\ degree\ vectors=} 81}
{\tt {\small
tried 246, got 81 pure-Sylv formulae [m,dimK1,dimK0]:\\
}
> sort(convert(\%,list),sort\_fnc);\\
}
\begin{eqnarray*}
\lefteqn{[[8, \,5, \,3, \,1120, \,1080],
 \,[8, \,3, \,5, \,1120, \,1080],}\\
 & & [6, \,9, \,3, \,1200, \,1120],\, [6, \,3, \,9, \,1200, \,1120],\\
 & & [4, \,9, \,7, \,1440, \, 1200], \,[4, \,7, \,9, \,1440, \,1200],\\
 & & [8, \,6, \,3, \,1400, \,1260],
 \,[8, \,3, \,6, \,1400, \, 1260],\, \dots]\\
\end{eqnarray*}
{\tt
> allDetVecs(Ns,Ds):\\
> allmsrtd := sort(convert(\%,list),sort\_fnc);\\
}
\centerline{$\mathit{From}, \,[-4, \,-3, \,-3], \,\mathit{to}, \,
[9,\,10,\,10],\, \mathit{start\ at},\,[-4,\, -3,\,-3]$}
{\tt\small
Tested 1452 m-vectors: assuming Pk's nonempty.\\
Found 488 det'l m-vecs, listed with matrix dim:\\
}
\begin{eqnarray*}
\lefteqn{\mathit{allmsrtd} := [[6,3,1, \,224],
 \,[6,1,3, \,224], \,[1,7,5, \,224],}\\
&&[1,5,7, \,224], \,[3,7,1, \,240],\,
    [4,7,1, \,240], \,[3,1,7, \,240],\\
&&[4,1,7, \,240], \,[6,3,0, \,262], \,
    [6,0,3, \,262], \, [1,8,5, \,262],\\
&&[1,5,8, \,262], \,\dots]\mbox{\hspace{154pt}}\\
\end{eqnarray*}
{\tt
> for i from 1 to nops( allmsrtd ) do print (\\
allmsrtd[i],Pksets(Ns,Ds,vector([allmsrtd[i][1],\\
allmsrtd[i][2],allmsrtd[i][3]]))): od:
}
\begin{center}
$[6, \,3, \,1, \, 224], \,[[4, \,4], \,[2, \,2], \,[1, \,1]]$\\
$[ 6, \, 1, \, 3, \, 224 ] , \, [ [ 4, \,4], \,[1, \,1], \,[2, \,2]] $\\
$[1, \,7, \,5, \, 224], \, [[1, \,1], \,[4, \,4], \,[3, \,3]]$\\
$[1, \,5, \,7, \, 224], \, [[1, \,1], \,[3, \,3], \,[4, \,4]]$\\
$[3, \,7, \,1, \,240], \,[[2, \,2], \,[4, \,4], \,[1, \,1]]$\\
$[4, \,7, \,1, \,240], \,[[3, \,3], \,[4, \,4], \,[1, \,1]]$\\
$[3, \,1, \,7, \,240], \,[[2, \,2], \,[1, \,1], \,[4, \,4]]$\\
$[4, \,1, \,7, \,240], \,[[3, \,3], \,[1, \,1], \,[4, \,4]]$\\
$[6, \,3, \,0, \,262], \,[[4, \,4], \,[2, \,2], \,[0, \,\mathit{
NO\_INT}, \,0]]$\\
$\dots$
\end{center}
The last 2 commands find all 488 determinantal vectors.
The smallest formulae are indicated
(the minimum dimension is 224) and for some
we report the $P_k$'s.
No determinantal formulae is pure Sylvester.
Notice that the assumption of empty sets $P_k$ is used only in
order to bound the search, but within the appropriate range
empty $P_k$'s are considered, so no valid degree vector is missed.
This is illustrated by the last $P_3=\emptyset$ marked {\em NO\_INT}.

The vectors predicted by Theorem \ref{Tdetmvec} are found among
those produced above, including
$m^{(123)}=(6,5,3),\, m^{(213)}=(4,9,3),\, m^{(312)}=(0,9,7)$.
The corresponding matrix dimensions are 672, 600, and 800.
\end{oneshot}

\bigskip

\begin{oneshot}{Example \protect\ref{Ex12_21} (Cont'd)}
Recall that  $l=(1,2), d=(2,1)$.
The only pure B\'ezout formulae are the
2 determinantal formulae of Ex.\ \ref{Ex12_21},
for which we have $m^\pi=(5,0), m^{\pi'}=(1,1)$.\\
{\tt
> Ns:=vector([1,2]):Ds:=vector([2,1]):\\
> summs:=allsums(Ns):\\
> findBez(Ns,Ds,true); \#not only determinantal\\
}
\centerline{$\mathit{low-upper\ bounds,\ 1st\ candidate:},
\,[0, \,0], \,[6, \,1], \, [0, \,0]$}
{\tt {\small
Searched degree m-vecs for ANY pure Bezout formula.\\
Tested 15, found 2 pure-Bezout [m,dimK0,dimK1]:\\
}}
\centerline{$\{[5, \,0, \,6, \,6], \,[1, \,1, \,6, \,6]\}$}
The search examined 15 degree vectors between the shown bounds.
It is clear that both vectors are determinantal because the matrix
dimensions are for both $6\times 6$.
\end{oneshot}

\section{Further work}

Our results can be generalized to polynomials with scaled
supports or with a different degree $d$ per polynomial.
One question is whether the vectors $m',m''$ of Definition
\ref{Dpertm} lead to smaller or larger matrices than $m$.
Notice that certain cohomologies, which were nonzero for $m$,
may vanish for $m'$ or $m''$.
We plan to complete the description of hybrid
determinantal formulae. We would also like to answer in general
the question stated in \S \ref{degreevectors} of determining a priori
the degree
vectors yielding the smallest determinantal formulae in all
possible cases.
A problem related to the Sylvester formulae calls for
identifying in advance the nonzero maximal minor in the matrix,
which leads to finding a determinant with exact degree in some
polynomial.

\end{document}